\documentclass[12pt]{article}
\usepackage{e-jc}
\usepackage{amsmath,tikz}
\usepackage{enumerate}

\newcommand{\F}{\mathcal{F}}
\newcommand{\R}{\mathcal{R}}
\newcommand{\N}{\mathbb{N}}
\newcommand{\Z}{\mathbb{Z}}
\newcommand{\into}{\hookrightarrow}
\newcommand{\pto}{\xrightarrow{poi}}
\newcommand{\Cl}{C^\ell}
\newcommand{\Clstar}{C^{\ell^*}}

\dateline{Nov 28, 2020}{Feb 8, 2021}{TBD}


\Copyright{The authors. Released under the CC BY-ND license (International 4.0).}

\title{On Ramsey-minimal infinite graphs}

\author{Jordan Mitchell Barrett\\
\small School of Mathematics and Statistics\\[-0.8ex]
\small Victoria University of Wellington\\[-0.8ex]
\small Wellington, New Zealand\\
\small\tt math@jmbarrett.nz\\
\and
Valentino Vito\\
\small Department of Mathematics\\[-0.8ex]
\small Universitas Indonesia\\[-0.8ex]
\small Depok, Indonesia\\
\small\tt valentino.vito@sci.ui.ac.id}

\begin{document}

\maketitle

\vspace{-2cm}

\begin{abstract}
For fixed finite graphs $G$, $H$, a common problem in Ramsey theory is to study graphs $F$ such that $F \to (G,H)$, i.e.\ every red-blue coloring of the edges of $F$ produces either a red $G$ or a blue $H$. We generalize this study to infinite graphs $G$, $H$; in particular, we want to determine if there is a minimal such $F$. This problem has strong connections to the study of \textit{self-embeddable} graphs: infinite graphs which properly contain a copy of themselves. We prove some compactness results relating this problem to the finite case, then give some general conditions for a pair $(G,H)$ to have a Ramsey-minimal graph. We use these to prove, for example, that if $G=S_\infty$ is an infinite star and $H=nK_2$, $n \ge 1$ is a matching, then the pair $(S_\infty,nK_2)$ admits no Ramsey-minimal graphs.
\end{abstract}

\maketitle

\section{Introduction}

Let $F$, $G$ and $H$ be possibly infinite, simple graphs with no isolated vertices. We follow some notation in \cite{schaefer}. We say that $F$ \textit{arrows} $(G,H)$ or that $F \to (G,H)$ if for every red-blue coloring of the edges of $F$, there exists either a red $G$ or a blue $H$ contained in $F$. In this case, we say that $F$ is an \emph{$(G,H)$-arrowing graph}. A red-blue coloring of $F$ is called \emph{$(G,H)$-good} if $F$ does not contain a red $G$ or a blue $H$ with respect to the coloring. An alternate definition for $F \to (G,H)$ would then be that the graph $F$ admits no $(G,H)$-good coloring.

A $(G,H)$-arrowing graph $F$ is said to be \emph{$(G,H)$-minimal} if there is no proper subgraph $F' \subset F$ such that $F' \to (G,H)$. In other words, $F$ is $(G,H)$-minimal if it arrows $(G,H)$ and $F-e \not\to (G,H)$ for every $e \in E(F)$. The collection of all $(G,H)$-minimal graphs is denoted as $\R(G,H)$, and it satisfies the symmetric property $\R(G,H)=\R(H,G)$.

The problem involving $(G,H)$-minimal graphs is classically done for finite $G$ and $H$, as introduced in \cite{burr-type}. One of the major problems that arose was determining whether $\R(G,H)$ is finite or infinite. Following the studies done by Beardon \cite{beardon} on magic labelings, C\'aceres \textit{et al.} \cite{caceres} on metric dimensions, and Stein \cite{stein} on extremal graph theory, we attempt to extend this finite problem to an infinite one. To our knowledge, this is the first serious attempt to do so. It appears that some properties which are expected to be true for finite graphs do not hold in the scope of infinite graphs.

For finite graphs $G$ and $H$, it is known that $\R(G,H)$ is nonempty. This is because we can obtain a $(G,H)$-minimal graph from an arbitrary $(G,H)$-arrowing graph by iteratively deleting enough edges. However, if one of $G$ or $H$ is infinite, then $\R(G,H)$ might be empty. As we shall see in Example \ref{exm-empty-ramsey}, the ray $P_\infty$ and $K_2$ as a pair do not admit any minimal graph. If we consider the double ray $P_{2\infty}$ instead of $P_\infty$, we have that $\R(P_{2\infty},K_2)=\{P_{2\infty}\}$, and thus a minimal graph exists. An intriguing but difficult problem in general would be to classify which pairs $(G,H)$ induce an empty (resp., nonempty) $\R(G,H)$.

\begin{problem}\label{main-problem}
For which pairs of graphs $(G,H)$ is $\R(G,H)$ empty?
\end{problem}

The study of Ramsey-minimal properties of infinite graphs is naturally related to graphs which are isomorphic to some proper subgraph of themselves. We will call such graphs \emph{self-embeddable}. Note that if $F$ is a self-embeddable graph, then we can pick an $F' \subset F$ isomorphic to $F$. Thus, if $F \to (G,H)$ is self-embeddable, then it is not $(G,H)$-minimal since we can choose a proper subgraph $F' \to (G,H)$.

\begin{observation}\label{main-obs}
A $(G,H)$-arrowing graph $F$ that is self-embeddable cannot be minimal.
\end{observation}

The notion of a self-embeddable graph differs from that of a \emph{self-contained graph} \cite{shekarriz}, which is one isomorphic to a proper \emph{induced} subgraph of itself. While self-contained graphs have applications in other problems, such as the tree alternative conjecture \cite{bonato}, we do not require the proper subgraph to be induced in our case, hence the differing vocabulary.

The general outline of this paper is as follows. In Section \ref{prelim}, we present the common notation and conventions that we use for this paper. We then give some compactness results for Ramsey-minimal graphs in Section \ref{compactness}. In Section \ref{general}, we obtain some general progress on Problem \ref{main-problem}. In Section \ref{matchings}, we turn our attention to the case where $G$ is an infinite graph and $H$ is a \emph{matching} $nK_2$. For an example of previous work on Ramsey-minimal finite graphs involving matchings, see Burr \textit{et al.} \cite{burr-matchings}.

\section{Preliminaries}\label{prelim}

In this paper, we exclusively work with simple graphs $G=(V(G),E(G))$ with no isolated vertices (i.e.\ every vertex of $G$ is adjacent to another vertex). Our graphs are taken to be countable (including finite), with the exception of the graphs of Section \ref{compactness} which may be uncountable. Let $\N=\{1,2,\ldots\}$ be the set of natural numbers. For $n \in \N$, $nG$ denotes the graph consisting of $n$ disjoint copies of $G$.

We say that $H$ is a \emph{subgraph} of $G$ (or simply $H \subseteq G$) if $V(H) \subseteq V(G)$ and $E(H) \subseteq E(G)$. A subgraph $H$ of $G$ is \emph{proper}, and written as $H \subset G$, if $E(H)$ is a proper subset of $E(G)$. Also, we say that $G$ (\emph{properly}) \emph{contains} $H$, or that $H$ (\emph{properly}) \emph{embeds} into $G$, if there is a (proper) subgraph $H'$ of $G$ such that $H' \cong H$.

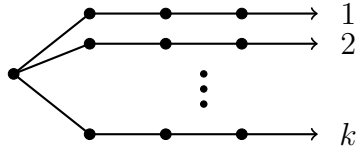
\begin{figure}
    \centering
    \begin{tikzpicture}[x=10mm,y=4mm]
        \draw[fill=black] (0,0) circle (2pt);
        \draw[fill=black] (1,1) circle (2pt);
        \draw[fill=black] (2,1) circle (2pt);
        \draw[fill=black] (3,1) circle (2pt);
        \draw[fill=black] (1,2) circle (2pt);
        \draw[fill=black] (2,2) circle (2pt);
        \draw[fill=black] (3,2) circle (2pt);
        \draw[fill=black] (1,-2) circle (2pt);
        \draw[fill=black] (2,-2) circle (2pt);
        \draw[fill=black] (3,-2) circle (2pt);
        
        \draw[thick] (0,0) -- (1,1) -- (2,1) -- (3,1);
        \draw[thick] (0,0) -- (1,2) -- (2,2) -- (3,2);
        \draw[thick] (0,0) -- (1,-2) -- (2,-2) -- (3,-2);
        \draw[->,thick] (3,1) -- (4,1);
        \draw[->,thick] (3,2) -- (4,2);
        \draw[->,thick] (3,-2) -- (4,-2);
        
        \draw[fill=black] (2.5,0) circle (1.2pt);
        \draw[fill=black] (2.5,-0.5) circle (1.2pt);
        \draw[fill=black] (2.5,-1) circle (1.2pt);
        
        \node at (4.4,2) {$1$};
        \node at (4.4,1) {$2$};
        \node at (4.4,-2) {$k$};
    \end{tikzpicture}
    
    \caption{A $k$-ray $P_{k\infty}$.}
    \label{fig:P_kinf}
\end{figure}

The \emph{ray} $P_\infty$ is an infinite graph of the form $(\{x_0,x_1,\ldots\},\{x_0x_1,x_1x_2,\ldots\})$, where $x_0$ is its \textit{endpoint}. The \emph{double ray} $P_{2\infty}$, on the other hand, is of the form $(\{x_n: n \in \Z\},\{x_nx_{n+1}:n \in \Z\})$. In general, a \emph{$k$-ray} $P_{k\infty}$ (shown in Figure \ref{fig:P_kinf}), $k \ge 1$, is formed by identifying the endpoints of $k$ distinct rays.

A family of graphs of particular interest is the family of comb graphs. Let $\ell\colon \N \to \N$ be a function. The \emph{comb} $\Cl$ is a graph obtained from a base ray $P_\infty$ (called the \emph{spine}) by attaching, for every $n$, a path $P_{\ell(n)}$ of order $\ell(n)$ by one of its endpoints to the vertex $x_n$ of $P_\infty$. Other infinite graphs of interest include the (countably) infinite complete graph $K_\infty$ and the (countably) infinite star $S_\infty=K_{1,\infty}$.

We use some terminologies of embeddings from \cite{binns,halin,hamann}. Recall that a \emph{graph homomorphism} $G \to H$ is a map $\varphi\colon V(G) \to V(H)$ such that if $vw \in E(G)$, then $\varphi(v)\varphi(w) \in E(H)$. If $\varphi$ is injective, then the homomorphism is called an \textit{embedding} $G \into H$. An embedding $G \into G$ is said to be a \emph{self-embedding} of $G$. A self-embedding is \emph{nontrivial} if its image, seen as a graph with vertex set $\varphi(V(G))$ and edge set $\{\varphi(v)\varphi(w):vw \in E(G)\}$, is a proper subgraph of $G$.

A graph $G$ is said to be \emph{self-embeddable} if it has a nontrivial self-embedding. In other words, a self-embeddable graph is a graph that properly embeds into itself. We say that $G$ is \emph{strongly self-embeddable} if it admits an embedding into $G-v$ for every $v \in V(G)$. A strongly self-embeddable graph is clearly self-embeddable, but the converse does not hold in general, as shown in the following example.

\begin{figure}[b]
    \centering
    \begin{tikzpicture}[x=10mm,y=10mm]
        \draw[fill=red] (0,0) circle (2pt);
        \draw[fill=red] (1,0) circle (2pt);
        \draw[fill=red] (2,0) circle (2pt);
        \draw[fill=red] (3,0) circle (2pt);
        \draw[fill=black] (0,1) circle (2pt);
        \draw[fill=red] (1,1) circle (2pt);
        \draw[fill=red] (2,1) circle (2pt);
        \draw[fill=red] (3,1) circle (2pt);
        \draw[fill=red] (-1,0) circle (2pt);
        \draw[fill=red] (-2,0) circle (2pt);
        \draw[fill=red] (-3,0) circle (2pt);
        
        \draw[thick,red] (-3,0) -- (-2,0) -- (-1,0) -- (0,0) -- (1,0) -- (2,0) -- (3,0);
        \draw[thick] (0,0.08) -- (0,1);
        \draw[thick,red] (1,0) -- (1,1);
        \draw[thick,red] (2,0) -- (2,1);
        \draw[thick,red] (3,0) -- (3,1);
        \draw[->,thick,red] (-3,0) -- (-4,0);
        \draw[->,thick,red] (3,0) -- (4,0);
        
        \node at (0,-0.3) {$v$};
        
    \end{tikzpicture}
    
    \caption{A self-embeddable graph $G$ that is not strongly self-embeddable. The red subgraph illustrates the image of a nontrivial self-embedding of $G$ in the form of a right translation by $1$.}
    \label{fig:not-strongly-emb}
\end{figure}
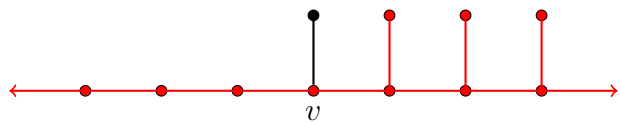

\begin{example}
The infinite star $S_\infty$ is self-embeddable but not strongly so (since $S_\infty$ does not embed into the null graph $S_\infty-c$, where $c$ is its center vertex). Another example can be found in the graph $G$ of Figure \ref{fig:not-strongly-emb}. It is self-embeddable, with a right translation as its nontrivial self-embedding. However, $G$ does not embed into the disconnected graph $G-v$, where $v$ is the vertex indicated in the figure.
\end{example}

By induction, we easily obtain the following stronger properties for strongly self-embeddable graphs.

\begin{proposition}\label{first-prop}
Let $G$ be strongly self-embeddable. Then:
\begin{enumerate}[\normalfont(i)]
\item $G-V$ contains $G$ for every finite $V \subset V(G)$;
\item $G-E$ contains $G$ for every finite $E \subset E(G)$.
\end{enumerate}
\end{proposition}

\begin{proof}
(i) Suppose that $V=\{v_1,\ldots,v_n\}$ and that $H=G-\{v_1,\ldots,v_{n-1}\}$ contains $G$. In other words, there is a subgraph $G' \subseteq H$ isomorphic to $G$. Since $G'$ and $G$ are isomorphic, $G'$ is strongly self-embeddable. Thus, $G'-v_n$ must contain $G' \cong G$. Since $G'-v_n \subseteq H-v_n = G-V$, we can conclude that $G-v$ contains $G$.

(ii) Suppose that $e=uv$. Since $G-v$ contains a copy of $G$ by hypothesis, we have that $G-e \supseteq G-v$ contains $G$ as well. The statement for any $G-E$ then follows by induction.
\end{proof}

\section{Compactness results}\label{compactness}

The aim of this section is, for infinite (possibly uncountable) graphs $F$, $G$, $H$, to express what $F \to (G,H)$ means in terms of their finite subgraphs $\hat{F}$, $\hat{G}$, $\hat{H}$.

Theorem \ref{thm:minim-finite} confirms that all $(G,H)$-minimal graphs are finite if $G$ and $H$ are finite. This result assures us that we only need to deal with the case where one of $G$ and $H$ is infinite, since taking both as finite graphs would produce a completely finite problem. We prove the theorem by topological means using Tychonoff's theorem.

Recall that, for a family of topological spaces $S_i$, $i \in I$, the \textit{product topology} on $\prod S_i$ is generated by basic open sets of the form $\prod U_i$, where each $U_i$ is open in $S_i$, and $U_i=S_i$ except for finitely many values of $i$. Tychonoff's theorem states that whenever each $S_i$ is compact, $\prod S_i$ is also compact.

\begin{theorem}\label{thm:minim-finite}
Let $F$ be a graph, and $G$ and $H$ be finite graphs. If $F \to (G,H)$, then there is a finite $\hat{F} \subseteq F$ such that $\hat{F} \to (G,H)$.
\end{theorem}

\begin{proof}
Let $X$ be the product of $|E(F)|$-many copies of the discrete space $S = \{ \text{red},\text{blue} \}$, equipped with the product topology. We can identify $X$ with the set of all functions $E(F) \to S$, i.e.\ the set of all red-blue colorings of $F$'s edges. By Tychonoff's theorem, $X$ is compact.

By assumption, $F \to (G,H)$, so for every coloring $c \in X$, we can pick a finite set of edges $D_c$ forming either a red $G$ or a blue $H$. Then, each $c|_{D_c}$ determines a basic open set $$O_c:=\{d\colon E(F) \to S: d|_{D_c}=c|_{D_c}\}.$$ Since $c \in O_c$ for all $c \in X$, the collection $\{ O_c: c \in X \}$ covers $X$. By compactness, there is a finite sequence $c_1, \ldots, c_n \in X$ so that $O_{c_1} \cup \cdots \cup O_{c_n} = X$.

Let $\hat{F}$ be the subgraph of $F$ induced by $D_{c_1} \cup \cdots \cup D_{c_n}$. We claim $\hat{F} \to (G,H)$. Pick a coloring $\hat{d}\colon E(\hat{F}) \to S$, and extend it arbitrarily to a coloring $d\colon E(F) \to S$. Then, there is an $i \leq n$ such that $d \in O_{c_i}$, so $D_{c_i} \subseteq E(\hat{F})$ either forms a red $G$ or a blue $H$ under $\hat{d}$, as required.
\end{proof}

Now, we want to be able to characterize embeddability of a graph $G$ into another graph $F$ in terms of embeddability of finite subgraphs $\hat{G} \subseteq G$ into $F$. We might want to prove something such as:
\begin{center}
    $G$ embeds into $F$ if and only if every finite subgraph $\hat{G} \subseteq G$ embeds into $F$.
\end{center}

\begin{figure}
    \centering
    \begin{tikzpicture}[x=5mm,y=5mm]
        \draw[fill=black] (0,0) circle (2pt);
        \draw[fill=black] (1,-1) circle (2pt);
        \draw[fill=black] (1,1) circle (2pt);
        \draw[fill=black] (2,0) circle (2pt);
        \draw[fill=black] (3,-1) circle (2pt);
        \draw[fill=black] (3,1) circle (2pt);
        \draw[fill=black] (4,0) circle (2pt);
        \draw[fill=black] (5,-1) circle (2pt);
        \draw[fill=black] (5,1) circle (2pt);
        \draw[fill=black] (6,0) circle (2pt);
        \draw[fill=black] (7,-1) circle (2pt);
        \draw[fill=black] (7,1) circle (2pt);
        \draw[fill=black] (8,0) circle (2pt);
        
        \draw[thick] (8.7,0.7) -- (7,-1) -- (5,1) -- (3,-1) -- (1,1) -- (0,0) -- (1,-1) -- (3,1) -- (5,-1) -- (7,1) -- (8.7,-0.7);
        
        \draw[fill=black] (9.5,0) circle (1.2pt);
        \draw[fill=black] (10,0) circle (1.2pt);
        \draw[fill=black] (10.5,0) circle (1.2pt);
        
        \node at (5.5,-2) {$G$};

        \draw[fill=black] (-5,0) circle (2pt);
        \draw[fill=black] (-4,-1) circle (2pt);
        \draw[fill=black] (-4,1) circle (2pt);
        \draw[fill=black] (-3,0) circle (2pt);
        
        \draw[fill=black] (-5,-2) circle (2pt);
        
        \draw[fill=black] (-5,-4) circle (2pt);
        \draw[fill=black] (-4,-5) circle (2pt);
        \draw[fill=black] (-4,-3) circle (2pt);
        \draw[fill=black] (-3,-4) circle (2pt);
        \draw[fill=black] (-2,-5) circle (2pt);
        \draw[fill=black] (-2,-3) circle (2pt);
        \draw[fill=black] (-1,-4) circle (2pt);
        
        \draw[fill=black] (-5,-6) circle (2pt);
        
        \draw[fill=black] (-5,-8) circle (2pt);
        \draw[fill=black] (-4,-9) circle (2pt);
        \draw[fill=black] (-4,-7) circle (2pt);
        \draw[fill=black] (-3,-8) circle (2pt);
        \draw[fill=black] (-2,-9) circle (2pt);
        \draw[fill=black] (-2,-7) circle (2pt);
        \draw[fill=black] (-1,-8) circle (2pt);
        \draw[fill=black] (0,-9) circle (2pt);
        \draw[fill=black] (0,-7) circle (2pt);
        \draw[fill=black] (1,-8) circle (2pt);
        
        \draw[thick] (-5,0) -- (-4,-1) -- (-3,0) -- (-4,1) -- (-5,0) -- (-5,-2) -- (-5,-4) -- (-4,-5) -- (-3,-4) -- (-2,-5) -- (-1,-4) -- (-2,-3) -- (-3,-4) -- (-4,-3) -- (-5,-4) -- (-5,-6) -- (-5,-8) -- (-4,-9) -- (-3,-8) -- (-2,-9) -- (-1,-8) -- (0,-9) -- (1,-8) -- (0,-7) -- (-1,-8) -- (-2,-7) -- (-3,-8) -- (-4,-7) -- (-5,-8) -- (-5,-10);
        
        \draw[fill=black] (-3,-10) circle (1.2pt);
        \draw[fill=black] (-3,-10.5) circle (1.2pt);
        \draw[fill=black] (-3,-11) circle (1.2pt);
        
        \node at (0,-5.5) {$F$};
    \end{tikzpicture}
    
    \caption{Every finite subgraph $\hat{G} \subseteq G$ embeds into $F$, but $G$ itself does not.}
    \label{fig:cpct-embed}
\end{figure}
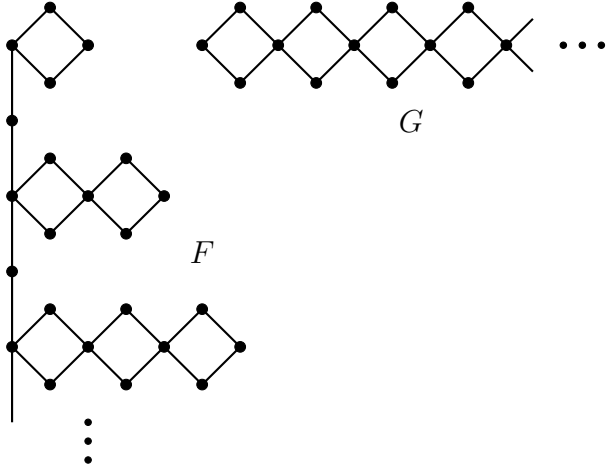

However, this statement is not true; a counterexample is shown in Figure \ref{fig:cpct-embed}. The problem is that the embeddings $\hat{G} \into F$ are incompatible in some sense---larger finite subgraphs $\hat{G} \subseteq G$ must be embedded further down $F$, and so there is no way to ``stitch together'' these embeddings to get an embedding $G \into F$. To ensure compatibility between partial embeddings, we instead work with the following notion of a \emph{pointed graph}. We also need to ensure that our graphs are \textit{locally finite}: that is, $\deg(v)<\infty$ for each vertex $v$.

A \emph{pointed graph} is a triple $G = (V,E,*)$, where $(V,E)$ is a graph, and $* \in V$ is a specified vertex of $G$, called the \emph{basepoint} of $G$. A \emph{pointed subgraph} $H$ is a subgraph of $G$ such that $*_H=*_G$. A \emph{pointed homomorphism} is a graph homomorphism mapping basepoints to basepoints---we call it a \emph{pointed embedding} if it is injective. A pointed graph is locally finite, connected, etc. if the underlying graph is.

Using pointed graphs, we can obtain a version of our desired result. We first provide the following form of K\H{o}nig's infinity lemma to pave way for the compactness argument used in the proof of Proposition \ref{prop:cpct-embed}.

\begin{lemma}[\protect{\cite[Lemma 8.1.2]{diestel}}]\label{konig}
Let $V_0, V_1, \ldots$ be an infinite sequence of disjoint nonempty finite sets, and let $K$ be a graph on their union. Assume that every vertex in $V_{n+1}$ has a neighbor in $V_n$. Then, $K$ contains a ray $v_0v_1\ldots$ such that $v_n \in V_n$ for all $n$.
\end{lemma}

\begin{proposition}\label{prop:cpct-embed}
Let $F$ and $G$ be locally finite, pointed graphs, and suppose $G$ is connected. If every connected, finite, pointed subgraph $\hat{G}$ of $G$ admits a pointed embedding into $F$, then $G$ admits a pointed embedding into $F$.
\end{proposition}

\begin{proof}
Since $G$ is locally finite and connected, it must be countable (see Exercise 1, Chapter 8 of \cite{diestel}). The lemma clearly holds if $G$ is finite, so assume that $G$ is countably infinite. Enumerate $V(G)$ as follows: let $v_0 = *$, let $v_1, \ldots, v_k$ be the neighbors of $*$, let $v_{k+1}, \ldots, v_\ell$ be the vertices of distance 2 to $*$, etc. This indeed enumerates $G$ by the assumption that $G$ is locally finite and connected. Then, the induced subgraphs $\hat{G}_n:=G[v_0, \ldots, v_n]$ are all connected, finite, pointed subgraphs of $G$.

For each $n$, let $V_n$ be the set of pointed embeddings $\hat{G}_n \into F$. Each $V_n$ is nonempty by assumption. Inductively, we show each $V_n$ is finite. We see that $V_0$ is finite, since there is a unique embedding $\hat{G}_0 \into F$ mapping $*_G$ to $*_F$. Now assume $V_n$ is finite, and pick $f \in V_n$. Since $\hat{G}_{n+1}$ is connected, $v_{n+1}$ is adjacent to some $v_j$ for $j \le n$. Since $F$ is locally finite, $f(v_j)$ has finite degree, so there are only finitely many ways to define $f(v_{n+1})$ and extend $f$ to an embedding in $V_{n+1}$. Since there are also only finitely many ways to choose $f \in V_n$, it follows that $V_{n+1}$ is finite, and so all the $V_n$ are finite by induction.

Similarly to before, let $K$ be the graph on $\bigcup_{n=0}^\infty V_n$, where we insert all edges between $f \in V_{n+1}$ and $f|_{\hat{G}_n} \in V_n$. By Lemma \ref{konig}, $K$ contains an infinite ray $f_0f_1\ldots$ such that $f_n \in V_n$ for all $n$. Define $f\colon V(G) \to V(F)$ by $f(v_n) = f_n(v_n)$. We claim $f$ is a pointed embedding $G \into F$.
\begin{enumerate}[(i)]
    \item $f$ is a graph homomorphism: suppose $v_nv_m \in E(G)$, where $n \le m$. We have $f(v_n) = f_n(v_n) = f_m(v_n)$ and $f(v_m) = f_m(v_m)$. Since $f_m$ is a graph homomorphism, $f(v_n)f(v_m)=f_m(v_n)f_m(v_m) \in E(F)$.
    
    \item $f$ is injective: suppose $f(v_n) = f(v_m)$ for $n \le m$. Then, $f_m(v_n) = f_n(v_n) = f_m(v_m) \implies v_n = v_m$ since $f_m$ is injective.
    
    \item $f$ is pointed: $f(*_G) = f_0(*_G) = *_F$ since $f_0$ is pointed.\qedhere
\end{enumerate}
\end{proof}

Interestingly, Proposition \ref{prop:cpct-embed} actually generalizes K\H{o}nig's lemma, in its more standard, graph-theoretical form (as found in \cite[Proposition 8.2.1]{diestel}).

\begin{corollary}[K\H{o}nig's lemma]
Every locally finite, connected, infinite graph $F$ contains a ray.
\end{corollary}

\begin{proof}
Pick an arbitrary basepoint $* \in F$. For every $n$, we claim that $P_n$ (with $*_{P_n}$ chosen as an endpoint) admits a pointed embedding into $F$. If $P_n$ does not admit a pointed embedding into $F$, then there is no vertex $v \in V(F)$ such that $d(*_F,v) \ge n$. Since $F$ is connected and locally finite, we can enumerate $V(F)$ as follows: let $v_0 = *$, let $v_1, \ldots, v_k$ be the neighbors of $*$, let $v_{k+1}, \ldots, v_\ell$ be the vertices of distance 2 to $*$, etc. By stage $n$, we will have enumerated all of $F$, hence $F$ is finite; contradiction. The result now follows from Proposition \ref{prop:cpct-embed}, with $G = P_\infty$ and $*_{P_\infty}$ as its endpoint.
\end{proof}

For pointed graphs $F$, $G$, and a non-pointed graph $H$, we write $F \pto (G,H)$ if for every red-blue coloring of the edges of $F$, there exists either a red $G$ as a pointed subgraph of $F$ or a blue $H$ in the underlying graph of $F$. This definition gives a stronger condition for $F$ than $F \to (G,H)$, and, unlike regular arrowing, $F \pto (G,H)$ and $F \pto (H,G)$ are not necessarily equivalent. We also note that $F\pto (G,H)$ only if $H \subseteq K_{1,\deg(*_F)}$.

\begin{theorem}\label{thm:pointed-arrowing}
Let $F$, $G$ be locally finite, pointed graphs and let $H$ be a graph. Suppose $G$ is connected and for every connected, finite, pointed $\hat{G} \subseteq G$, we have $F \pto (\hat{G},H)$. Then, $F \pto (G,H)$.
\end{theorem}

\begin{proof}
Take an arbitrary red-blue coloring of $F$ such that $F$ does not contain a blue $H$. Denote $F'$ as the pointed subgraph of $F$ induced by all the red edges. By assumption, all connected, finite, pointed $\hat{G} \subseteq G$ admits a pointed embedding into $F'$. By Proposition \ref{prop:cpct-embed}, $G$ admits a pointed embedding into $F'$, so a red $G$ exists as a pointed subgraph of $F$.
\end{proof}

Proposition \ref{prop:cpct-embed} is precisely Theorem \ref{thm:pointed-arrowing} for the case $H=K_2$. Thus, Theorem \ref{thm:pointed-arrowing} generalizes both Proposition \ref{prop:cpct-embed} and K\H{o}nig's lemma.

\section{General progress on Problem \ref{main-problem}}\label{general}

Throughout this section and the next, we fix a pair of (potentially infinite) graphs $G$ and $H$. We provide a sufficient condition under which $\R(G,H)$ is empty by first finding some suitable family of graphs $\F$ for the pair $(G,H)$.

\begin{theorem}\label{gen-thm-1}
Suppose that $\F$ is a (possibly infinite) collection of graphs such that:
\begin{enumerate}[\normalfont(1)]
    \item $F \to (G,H)$ for every $F \in \F$;
    \item Every $(G,H)$-arrowing graph $\Gamma$ contains some graph $F \in \F$.
\end{enumerate}
We have the following:
\begin{enumerate}[\normalfont(i)]
    \item $\R(G,H) \subseteq \{F \in \F: \text{$F$ is not self-embeddable}\}$;
    \item If every $F \in \F$ is self-embeddable, then $\R(G,H)$ is empty.
\end{enumerate}
\end{theorem}

\begin{proof}
(i) Fix a $(G,H)$-minimal graph $\Gamma$. By condition (2), there is an $F \in \F$ such that $F$ is contained in $\Gamma$. Since $F \to (G,H)$ by condition (1), we must have $F=\Gamma$ as $\Gamma$ is $(G,H)$-minimal. Therefore, $\Gamma \in \F$. By Observation \ref{main-obs}, $\Gamma$ is not self-embeddable, and we are done.

(ii) follows directly from (i).
\end{proof}

Conditions (1) and (2) are not sufficient to ensure that $\R(G,H)$ and $\{F \in \F: \text{$F$ is not self-embeddable}\}$ coincide. For example, take $G=P_\infty$, $H=K_2$ and $\F=\{P_\infty,P_{2\infty}\}$. While conditions (1) and (2) hold, $\R(P_\infty,K_2)$ can be shown to be empty while $P_{2\infty}$ is not self-embeddable. Hence, $\R(G,H) \subset \{F \in \F: \text{$F$ is not self-embeddable}\}$. That said, we can create an extra condition to make both sets equal.

\begin{theorem}\label{gen-thm-2}
Let $\F$ be a collection of graphs such that conditions (1) and (2) of Theorem \ref{gen-thm-1} hold. Suppose that we also have the following condition:
\begin{enumerate}[\normalfont(1)]
    \setcounter{enumi}{2}
    \item $F_1$ and $F_2$ do not contain each other for every different $F_1,F_2 \in \F$.
\end{enumerate}
We have the following:
\begin{enumerate}[\normalfont(i)]
    \item $\R(G,H)=\{F \in \F: \text{$F$ is not self-embeddable}\}$;
    \item $\R(G,H)$ is empty if and only if every $F \in \F$ is self-embeddable.
\end{enumerate}
\end{theorem}

\begin{proof}
(i) We prove that $\{F \in \F: \text{$F$ is not self-embeddable}\} \subseteq \R(G,H)$. Suppose $F \in \F$ is not self-embeddable. Since $F \to (G,H)$ by condition (1), it remains to show that no proper $F' \subset F$ arrows $(G,H)$. If we assume the contrary, then we have an $F'' \in \F$ contained in $F'$ by condition (2). This implies that $F$ contains $F''$ and contradicts condition (3). 

(ii) follows easily from (i).
\end{proof}

The preceding theorems have a few applications. For example, in order to prove Theorem \ref{star-vs-matching} of the next section, we will need to use Theorem \ref{gen-thm-1}(ii). Also, we can consider the special case where $\F$ is chosen as $\{G\}$. This yields the following results:

\begin{theorem}\label{gen-thm-3}
We have the following:
\begin{enumerate}[\normalfont(i)]
    \item $\R(G,K_2)$ is empty if and only if $G$ is self-embeddable. If it is nonempty, then $\R(G,K_2)=\{G\}$;
    \item If $H \neq K_2$, then $G \to (G,H)$ implies that $\R(G,H)$ is empty.
\end{enumerate}
\end{theorem}

\begin{proof}
(i) Take $\F=\{G\}$. Conditions (1)--(3) of Theorems \ref{gen-thm-1} and \ref{gen-thm-2} all hold when $H=K_2$. The statement directly follows from Theorem \ref{gen-thm-2}.

(ii) Again, take $\F=\{G\}$. Conditions (1) and (2) of Theorem \ref{gen-thm-1} are both satisfied. Now we just need to show that $G$ is self-embeddable. Color an arbitrary edge of $G$ blue and the rest of the edges red. Since $H \neq K_2$, there is no blue $H$ in $G$. So by the fact that $G \to (G,H)$, there is a red copy of $G$ in $G$. This proves that $G$ properly contains itself, and thus self-embeddable.
\end{proof}

The following examples demonstrate some direct applications of Theorem \ref{gen-thm-3} for some pairs of graphs:

\begin{example}\label{exm-empty-ramsey}
$\R(P_{k\infty},K_2)$ is empty if and only if $k=1$. When $k>1$, we have $\R(P_{k\infty},K_2)=\{P_{k\infty}\}$. These observations are obtained directly from Theorem \ref{gen-thm-3}(i).
\end{example}

\begin{example}
By Theorem \ref{gen-thm-3}(ii), $\R(S_\infty,K_{1,n})$ is empty for all $2 \le n \le \infty$ since $S_\infty \to (S_\infty,K_{1,n})$.
\end{example}

\begin{example}
$\R(K_\infty,H)$ is empty for all graphs $H$. This follows from Theorem \ref{gen-thm-3} since $K_\infty$ is self-embeddable and $K_\infty \to (K_\infty,H)$ for all $H$ by the infinite Ramsey theorem.
\end{example}

\section{Some results involving matchings}\label{matchings}

We saw in Theorem \ref{gen-thm-3}(i) an answer to Problem \ref{main-problem} whenever $H=K_2$. Now, let us consider the more general case where $H=nK_2$. It becomes apparent that the characteristics of $\R(G,nK_2)$, $n \ge 2$, are still related to whether $G$ is (strongly) self-embeddable.

\begin{theorem}\label{matchings-main}
We have the following:
\begin{enumerate}[\normalfont(i)]
    \item If $G$ is connected and not self-embeddable, then for all $n \ge 2$, we have $nG \in \R(G,nK_2)$ so that $\R(G,nK_2)$ is nonempty;
    \item If $G$ is strongly self-embeddable, then $\R(G,nK_2)$ is empty for all $n \ge 2$.
\end{enumerate}
\end{theorem}

\begin{proof}
(i) Fix a red-blue coloring of $nG$ which does not create a blue $nK_2$. It follows that there must be a component of $nG$ isomorphic to $G$ which is colored all red. Hence, $nG \to (G,nK_2)$.

Now let $e \in E(nG)$ be an arbitrary edge located in some component $G'$ of $nG$. We show that $nG-e \not\to (G,nK_2)$ by constructing a $(G,nK_2)$-good coloring of $nG-e$ as follows: for every component of $nG$ other than $G'$, color one of its edges blue; color the rest of the edges red. Since $G$ is connected, there are exactly $n-1$ components in $nG$ other than $G'$, so this coloring only manages to produce a blue $(n-1)K_2$. Also, since $G$ is not self-embeddable, there cannot be a red $G$ in any of the components of $nG$. By the connectivity of $G$, there cannot be a red $G$ in all of $nG$ either. Therefore, this coloring is indeed $(G,nK_2)$-good.

(ii) By appealing to Theorem \ref{gen-thm-3}(ii), it suffices to prove that $G \to (G,nK_2)$. Fix a red-blue coloring of $G$ which does not create a blue $nK_2$. We claim that this coloring creates a red $G$. We construct a set of vertices $V \subset V(G)$ using the following algorithm:
\begin{enumerate}[1.]
    \setlength{\itemsep}{0pt}
    \setlength{\parskip}{0pt}
    \item \textbf{initialize} $V=\emptyset$
    \item \textbf{while} $G-V$ contains a blue edge \textbf{do}
    \item \hspace{0.6cm} choose a blue edge $e=uv$ in $G-V$
    \item \hspace{0.6cm} $V \leftarrow V \cup \{u,v\}$
    \item \textbf{output} $V$
\end{enumerate}
This algorithm must terminate after at most $n-1$ while loop iterations since the edge chosen at each iteration must be independent from the edges chosen at previous iterations. It is then clear that the output $V$ is finite and that $G-V$ only contains red edges. By Proposition \ref{first-prop}(i), we have a copy of $G$ in $G-V$. It follows that there exists a red copy of $G$ in $G$, and we are done.
\end{proof}

\begin{remark}
We note that Theorem \ref{matchings-main}(i) can fail to hold if $G$ is disconnected. For example, given $G=2P_{2\infty}$, it can be shown that $(n+1)P_{2\infty} \to (G,nK_2)$. This implies that $(2n)P_{2\infty}$ is not minimal for $n \ge 2$, so $\R(2P_{2\infty},nK_2)$ cannot be shown to be nonempty using the previous line of reasoning.
\end{remark}

\begin{example}
Observe that $P_\infty$ is strongly self-embeddable, while $P_{k\infty}$ is connected and not self-embeddable for $k>1$. By Theorem \ref{matchings-main}, we have for every $n \ge 2$, $\R(P_{k\infty},nK_2)$ is empty if and only if $k=1$.
\end{example}

We note that the converse of Theorem \ref{matchings-main}(ii) does not necessarily hold. There indeed exists a graph $G=S_\infty$ not strongly self-embeddable such that $\R(G,nK_2)$ is empty for all $n$.

\begin{theorem}\label{star-vs-matching}
For all $n \ge 1$, $\R(S_\infty,nK_2)$ is empty.
\end{theorem}

We prove Theorem \ref{star-vs-matching} by first defining a collection of graphs $\F_n$ such that conditions (1) and (2) of Theorem \ref{gen-thm-1} hold.

\begin{lemma}\label{lemma-star}
For every $n \ge 1$, define $\F_n$ to be the collection of all graphs $F$ satisfying the following two conditions:
\begin{enumerate}[\normalfont(a)]
    \item $F$ contains exactly $n$ vertices of infinite degree forming the set $X=\{x_1,\ldots,x_n\}$;
    \item If $uv \in E(F)$, then at least one of $u$, $v$ is an element of $X$.
\end{enumerate}
Then, $\F_n$ satisfies conditions (1) and (2) of Theorem \ref{gen-thm-1}, with $G=S_\infty$ and $H=nK_2$.
\end{lemma}

\begin{proof}
To prove condition (1), we show that $F \to (S_\infty,nK_2)$ for all $F \in \F_n$ by induction on $n$. The base case $n=1$ follows from the fact that $\F_1=\{S_\infty\}$ and $S_\infty \to (S_\infty,K_2)$. Now assume that every graph in $\F_n$ arrows $(S_\infty,nK_2)$, and let $F \in \F_{n+1}$ be arbitrary. Suppose that $c$ is a red-blue coloring of $F$ which produces no red $S_\infty$.

Pick an arbitrary vertex $u$ adjacent to $x_{n+1} \in X$ such that $u \notin X$ and the edge $x_{n+1}u$ is colored blue. This can be done since otherwise, $x_{n+1}$ is incident to infinitely many red edges. Observe that $F':=F-\{x_{n+1},u\}$ is an element of $\F_n$, so it contains a blue $nK_2$ with respect to the coloring $c|_{F'}$. Since this blue $nK_2$ and the blue $x_{n+1}u$ form an independent set of edges, there exists a blue $(n+1)K_2$ in $F$ with respect to $c$. This proves that $F \to (S_\infty,(n+1)K_2)$ and completes the induction.

Now we prove condition (2). Suppose that $\Gamma$ arrows $(S_\infty,nK_2)$. We claim that $\Gamma$ contains at least $n$ vertices of infinite degree. Assume that $Y$, where $|Y|<n$, is the set of vertices of $\Gamma$ having an infinite degree. By coloring all edges incident to a vertex in $Y$ blue and the rest of the edges red, we obtain a $(S_\infty,nK_2)$-good coloring of $\Gamma$; contradiction. Now, we can take arbitrary vertices $x_1,\ldots,x_n$ of $\Gamma$ of infinite degree. The subgraph of $\Gamma$ induced by all edges that are incident to at least one of $x_1,\ldots,x_n$ is an element of $\F_n$, therefore condition (2) holds.
\end{proof}

\begin{proof}[Proof of Theorem \ref{star-vs-matching}]
Let $n \ge 1$. Take $\F_n$ as defined in Lemma \ref{lemma-star}. Invoking Theorem \ref{gen-thm-1}(ii), we need to show that every $F \in \F_n$ is self-embeddable.

Let $F \in \F_n$ be arbitrary. We aim to define a nontrivial self-embedding $\varphi$ of $F$. Denote $U$ as the complement of $X$ (i.e.\ $U:=V(F)\setminus X$), and for all $1 \le i \le n$, let $\rho_i\colon U \to \{0,1\}$ be such that $\rho_i(u)=1$ iff $u$ is adjacent to $x_i$. Define $\rho(u)=(\rho_1(u),\ldots,\rho_n(u))$ for every $u \in U$. Since the image of $\rho$ is finite, then by the infinite pigeonhole principle, there exists an $x \in \{0,1\}^n$ such that $\rho^{-1}(x)$ is an infinite set, say $\rho^{-1}(x)=\{u_1,u_2,\ldots\}$. Define $\varphi\colon V(F) \to V(F)$ as
$$\varphi(v)=\begin{cases}
u_{i+1}, & v=u_i,\ i \ge 1, \\ v, & \text{otherwise}.
\end{cases}$$

It is clear that $\varphi$ is injective, and it is nontrivial since $u_1 \notin \varphi(V(F))$. It remains to show that $\varphi$ is a graph homomorphism. Suppose that $uv \in E(F)$. By condition (b), we can assume without loss of generality that $v \in X$, say $v=x_j$ for some $1 \le j \le n$. Since $v \notin U$, we must have $v \notin \rho^{-1}(x)$, so $\varphi(v)=v$. If $u \notin \rho^{-1}(x)$, then $\varphi(u)\varphi(v) = uv \in E(F)$, and we are done. So assume that $u=u_i$ for some $i \ge 1$. Since $uv=u_ix_j$ is an edge, we have that $\rho_j(u_i)=1$. Thus, we also have $\rho_j(u_{i+1})=1$ since $\rho(u_{i+1})=x=\rho(u_i)$. It follows that
$$\varphi(u)\varphi(v)=\varphi(u_i)\varphi(x_j)=u_{i+1}x_j \in E(F).$$
The proof that $\varphi$ is a nontrivial self-embedding is then complete.
\end{proof}

Now, let $\Cl$ be a comb with ray $x_0x_1\ldots$ as its spine and a path $P_\ell(n)$ of order $\ell(n)$ attached to every $x_n$. Suppose that the value of $s_{\Cl}:=\min_{\ell(n)>1} n$, which is equal to the smallest natural number $n$ such that $\deg(x_n)=3$, exists (that is, $\Cl$ is not a ray). We can assume that our combs satisfy $s_{\Cl} \ge \ell(s_{\Cl})-1$ without loss of generality. Indeed, if $s=s_{\Cl}$ is such that $s<\ell(s)-1$, then we can define a function
\begin{equation}\label{l*}\ell^*(n)=\begin{cases}
1, & 1 \le n < \ell(s)-1 \\ s+1, & n=\ell(s)-1, \\ \ell(n-\ell(s)+s+1), & n \ge \ell(s),
\end{cases}\end{equation}
so that $\Clstar \cong \Cl$ and $s^*=s_{\Clstar}$ satisfies $s^* \ge \ell^*(s^*)-1$. Here, $\Clstar$ is basically the comb obtained from $\Cl$ by exchanging the positions of $x_0\ldots x_s$ and $P_{\ell(s)}$ in the comb.

The following theorem gives equivalent formulations to the statement that $\R(\Cl,nK_2)$ is empty for all $n \ge 2$. In particular, the theorem gives an answer for Problem \ref{main-problem} whenever $G$ is a comb and $H$ is a matching.

\begin{theorem}\label{thm:comb}
Let $\Cl$ be a comb which is not a ray with $x_0x_1\ldots$ as its spine. Suppose that the value of $s=\min_{\ell(n)>1} n$ is at least $\ell(s)-1$. The following statements are equivalent:
\begin{enumerate}[\normalfont 1.]
    \item $\R(\Cl,nK_2)$ is empty for all $n \ge 2$;
    \item $\R(\Cl,nK_2)$ is empty for some $n \ge 2$;
    \item $\Cl$ is self-embeddable;
    \item $\Cl$ is strongly self-embeddable;
    \item There exists a $p \geq 1$ such that $\ell(n) \leq \ell(n+p)$ for all $n \in \N$.
\end{enumerate}
\end{theorem}

\begin{proof}
Implication $(1 \to 2)$ is trivial, while implications $(2 \to 3)$ and $(4 \to 1)$ are both consequences of Theorem \ref{matchings-main}. So we just need to prove $(5 \to 4)$ and $(3 \to 5)$.

$5 \to 4$: Let $v \in V(\Cl)$. If $v=x_0$, then by positively translating $\Cl$ by $p$ (so that $x_n \mapsto x_{n+p}$ and $P_{\ell(n)}$ maps into $P_{\ell(n+p)}$), we obtain an embedding $\Cl \into \Cl-v$. Otherwise, there exists a $k \in \N$ such that $v$ is located in the path $P_{\ell(k)}$ attached to $x_k$. In this case, positively translate $\Cl$ by $ap$, where $a$ is chosen such that $ap>k$, to obtain the desired embedding into $\Cl-v$.

$3 \to 5$: Let $\varphi\colon V(\Cl) \to V(\Cl)$ be a nontrivial self-embedding given by the assumption that $\Cl$ is self-embeddable. We cannot have $\varphi(x_0)=x_0$, since that would mean $\varphi$ is the identity map, hence not nontrivial. Thus, $\varphi(x_0)$ is located in the path $P_{\ell(k)}$, which is attached to $x_k$, for some $k \in \N$. Suppose that $P_{\ell(k)}:=x_ky_1\ldots y_{\ell(k)-1}$.

If $\varphi(x_0)=x_k$, then $\varphi$ must be a positive translation by $k$. Hence, statement 5 holds by taking $p=k$ since $\varphi$ maps each $P_{\ell(n)}$ into $P_{\ell(n+k)}$. So assume that $\varphi(x_0)=y_j$ for some $1 \le j \le \ell(k)-1$. This has an implication that $\ell(k)>1$, and thus $s \le k$ since $s=\min_{\ell(n)>1} n$. Observe that we cannot have $j>s$ since that would mean $\varphi(x_s)=y_{j-s}$ has degree $2$, which is less than $\deg(x_s)=3$; contradiction. In summary, we have the inequality $j \le s \le k$.

We claim that $j<k$. Assume for the sake of contradiction that $j=s=k$. Since we have established that $j \le \ell(k)-1$ and $s \ge \ell(s)-1=\ell(k)-1$, we then have $\ell(k)-1=j=s=k$. It follows that $\varphi$ must be the map
$$\varphi(v)=\begin{cases}
y_{j-n}, & v=x_n,\ 0 \le n \le s-1, \\ x_{s-m}, & v=y_m,\ 1 \le m \le \ell(k)-1, \\ v, & \text{otherwise},
\end{cases}$$
which is not nontrivial; contradiction.

Now we can define $p=k-j \ge 1$. We prove that $\ell(n) \leq \ell(n+p)$ for all $n \ge s$ (the case where $n<s$ is trivial since then $\ell(n)=1$).

\textit{Case} 1. $j=s$. We have previously established the inequalities $j \le \ell(k)-1$ and $s \ge \ell(s)-1$. We thus have
$$\ell(s) \le s+1 = j+1 \le \ell(k) = \ell(s+p).$$
In addition, $\varphi$ necessarily maps each $P_{\ell(n)}$, $n>s$, into $P_{\ell(n+(k-s))}$. Hence, we also have $\ell(n) \leq \ell(n+p)$ for all $n>s$.

\textit{Case} 2. $j<s$. We can see that $\varphi$ necessarily maps each $P_{\ell(n)}$, $n \ge s$, into $P_{\ell(n+(k-j))}$. It follows that $\ell(n) \leq \ell(n+p)$ for all $n \ge s$.

In both cases, we see that statement 5 holds for the chosen $p=k-j$. This completes the proof.
\end{proof}

\begin{figure}
    \centering
    \begin{tikzpicture}[x=10mm,y=7mm]
        \draw[fill=black] (0,0) circle (2pt);
        \draw[fill=red] (1,0) circle (2pt);
        \draw[fill=red] (2,0) circle (2pt);
        \draw[fill=red] (3,0) circle (2pt);
        \draw[fill] (2,1) circle (2pt);
        \draw[fill=red] (3,1) circle (2pt);
        \draw[fill] (3,2) circle (2pt);
        \draw[fill=red] (4,0) circle (2pt);
        \draw[fill=red] (4,1) circle (2pt);
        \draw[fill=red] (4,2) circle (2pt);
        \draw[fill] (4,3) circle (2pt);
        
        \draw[thick] (0,0) -- (0.92,0);
        \draw[thick,red] (1,0) -- (2,0) -- (3,0) -- (4,0);
        \draw[thick] (2,0.08) -- (2,1);
        \draw[thick,red] (3,0) -- (3,1);
        \draw[thick] (3,1.08) -- (3,2);
        \draw[thick] (4,2.08) -- (4,3);
        \draw[thick,red] (4,0) -- (4,1) -- (4,2);
        \draw[->,thick,red] (4,0) -- (5,0);
        
        \node at (0,-0.6) {$x_0$};
        \node at (1,-0.6) {$x_1$};
        \node at (2,-0.6) {$x_2$};
        \node at (3,-0.6) {$x_3$};
        \node at (4,-0.6) {$x_4$};
        
        \draw[fill=red] (6,0) circle (2pt);
        \draw[fill=red] (7,0) circle (2pt);
        \draw[fill=red] (8,0) circle (2pt);
        \draw[fill=red] (9,0) circle (2pt);
        \draw[fill=black] (7,1) circle (2pt);
        \draw[fill=red] (8,1) circle (2pt);
        \draw[fill=red] (9,1) circle (2pt);
        \draw[fill=black] (7,2) circle (2pt);
        \draw[fill=red] (10,0) circle (2pt);
        \draw[fill=red] (10,1) circle (2pt);
        
        \draw[thick,red] (6,0) -- (7,0) -- (8,0) -- (9,0) -- (10,0);
        \draw[thick] (7,0.08) -- (7,1) -- (7,2);
        \draw[thick,red] (8,0) -- (8,1);
        \draw[thick,red] (9,0) -- (9,1);
        \draw[thick,red] (10,0) -- (10,1);
        \draw[->,thick,red] (10,0) -- (11,0);
        
        \node at (6,-0.6) {$x_0$};
        \node at (7,-0.6) {$x_1$};
        \node at (8,-0.6) {$x_2$};
        \node at (9,-0.6) {$x_3$};
        \node at (10,-0.6) {$x_4$};
        
        \node at (2.5,-1.5) {(a)};
        \node at (8.5,-1.5) {(b)};
        
    \end{tikzpicture}
    
    \caption{Two combs $\Cl$ such that (a) $\ell(n)=n$, and (b) $\ell(1)=3$ and $\ell(n)=2$ for $n>1$. For each of the two combs, the red subgraph illustrates the graph image of its nontrivial self-embedding.}
    \label{fig:combs}
\end{figure}
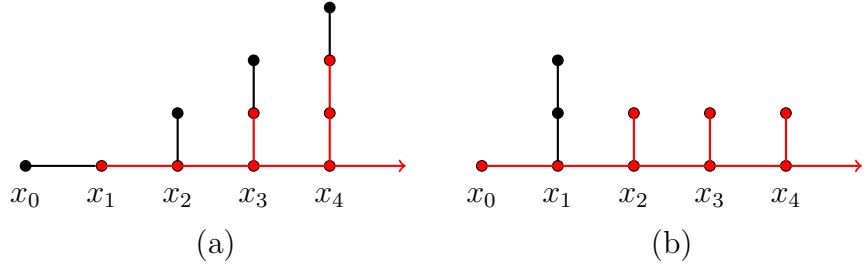

\begin{example}
If $\ell(n)=n$, then $\Cl$ satisfies $s \ge \ell(s)-1$ (with $s=2$) as well as statement 5 of Theorem \ref{thm:comb} for any choice of $p \ge 1$. As such, $\Cl$ is strongly self-embeddable via a positive translation. Figure \ref{fig:combs}(a) shows a positive translation of $\Cl$ by $1$. In addition, we have that $\R(\Cl,nK_2)$ is empty for all $n \ge 2$ by Theorem \ref{thm:comb}.
\end{example}

\begin{example}
Suppose that $\ell(1)=3$ and $\ell(n)=2$ for $n>1$. We have $s<\ell(s)-1$ (with $s=1$). So we define, using formula (\ref{l*}) preceding Theorem \ref{thm:comb}, a function $\ell^*$ such that
$$\ell^*(n)=\begin{cases}
1, & n=1 \\ 2, & n>1.
\end{cases}$$
The comb $\Clstar$ is illustrated in Figure \ref{fig:combs}(b) as a red subgraph of $\Cl$. We see that $\R(\Clstar,nK_2)$ is always empty since $\Clstar$ satisfies statement 5 of Theorem \ref{thm:comb}. By the fact that $\Cl \cong \Clstar$, we have that $\R(\Clstar,nK_2)$ is always empty as well.
\end{example}

\begin{example}
If $\ell(1)=2$ and $\ell(n)=1$ for $n>1$, then $s \ge \ell(s)-1$ (with $s=1$). However, statement 5 does not hold since $\ell(1)>\ell(1+p)$ for all $p \ge 1$. This implies that $\Cl$ is not self-embeddable and $\R(\Cl,nK_2)$ is nonempty for all $n \ge 2$. From Theorem \ref{matchings-main}(i), we can infer that $n\Cl \in \R(\Cl,nK_2)$ in this case.
\end{example}

\section{Concluding remarks}

Problem \ref{main-problem}, in its full generality, is quite a challenging problem to attack. For this reason, we chose to devote a significant part of this study to the particular case where $H$ is a matching. Even then, we were not able to completely answer Problem \ref{main-problem}. While Theorem \ref{matchings-main} managed to get us closer, we still have the following problem involving $\R(G,nK_2)$.

\begin{problem}
Is there a sufficient and necessary condition for $G$ under which $\R(G,nK_2)$ is empty for all $n \ge 2$?
\end{problem}

Further studies can also be done on other specific cases of the pair $(G,H)$. Of course, another avenue of research would be to consider multi-color Ramsey-minimal infinite graphs, as done in \cite{fox} for finite graphs.

While the compactness results of Section \ref{compactness} are interesting in their own right, our only application of them was to eliminate consideration of the case where both graphs $G$, $H$ are finite (Theorem \ref{thm:minim-finite}). We believe that these results, and other compactness results, could prove extremely useful in future studies of Ramsey-type problems for infinite graphs.

Apart from Section 3, we have only considered countable graphs in this paper. It would be interesting and worthwhile to study $\R(G,H)$ when at least one of $G$, $H$ is \textit{uncountable}. Presumably, this problem is significantly harder, and one would have to consider set-theoretic concerns.

\section*{Acknowledgements}

We would like to thank Fawwaz Fakhrurrozi Hadiputra for addressing Theorem \ref{gen-thm-3} in the beginning of our study, which inspired us to formulate Theorems \ref{gen-thm-1} and \ref{gen-thm-2}. Also, we would like to thank the anonymous referees for their valuable feedback and comments on this paper.

\end{document}